\documentclass[12pt]{amsart}
\usepackage{amscd,amssymb}
\usepackage[graph,frame,poly,arc]{xy}  
\usepackage[plainpages,backref,urlcolor=blue]{hyperref}

\theoremstyle{plain}
\newtheorem{thm}[subsection]{Theorem}
\newtheorem{lem}[subsection]{Lemma}
\newtheorem{prop}[subsection]{Proposition}
\newtheorem{cor}[subsection]{Corollary}

\theoremstyle{definition}
\newtheorem{rk}[subsection]{Remark}
\newtheorem{definition}[subsection]{Definition}
\newtheorem{ex}[subsection]{Example}

\newtheorem{question}[subsection]{Question}

\numberwithin{equation}{section}
\newcommand{\OO}{{\mathcal O}}

\newcommand{\F}{{\mathcal F}}

\newcommand{\C}{\mathbb{C}}
\newcommand{\PP}{\mathbb{P}}

\DeclareMathOperator{\rank}{rank}
\DeclareMathOperator{\im}{im}

\DeclareMathOperator{\codim}{codim}

\DeclareMathOperator{\depth}{depth }
\DeclareMathOperator{\syz}{syz}
\DeclareMathOperator{\Ext}{Ext}

\begin{document}

\title [Freeness   versus maximal degree of  singular subscheme for surfaces ]
{Freeness  versus maximal degree of the singular subscheme for surfaces in $\PP^3$ }

\author[Alexandru Dimca]{Alexandru Dimca$^1$}
\address{Univ. Nice Sophia Antipolis, CNRS,  LJAD, UMR 7351, 06100 Nice, France. }
\email{dimca@unice.fr}

\thanks{$^1$ Partially supported by Institut Universitaire de France.}



\begin{abstract} We show that a free surface in $\PP^3$ is characterized by the maximality of the degree of its singular subscheme, in the presence of an additional tameness condition. This is similar to the characterization of free plane curves by the maximality of their global Tjurina number given by A. A. du Plessis and C.T.C. Wall. Simple characterizations of  the nearly free tame surfaces are also given.

\end{abstract}
 
\maketitle


\section{Introduction} 

Let $D:f=0$ be a reduced hypersurface in the projective complex space $\PP^n$, defined by a homogeneous polynomial $f \in S=\C[x_0,...,x_n]$ of degree $d$. Let $\Sigma$ be the singular subscheme of $D$, defined by the Jacobian ideal $J_f$, which is the ideal in $S$ spanned by the first order partial derivatives $f_0,...,f_n$ of $f$ with respect to $x_0,...,x_n$. If $\dim \Sigma =m$, then the Hilbert polynomial $P(M(f))$ of the Milnor (a.k.a. Jacobian) algebra $M(f)=S/J_f$ has degree $m$, and the degree of $\Sigma$ is by definition $m! \cdot a$, where $a$ is the leading coefficient of $P(M(f))$. The minimal degree of a Jacobian relation for $f$ is the integer $mdr(f)$
defined to be the smallest integer $q\geq 0$ such that there is a nontrivial relation
\begin{equation}
\label{rel_m}
 \sum_{j=0,n}a_jf_j=0
\end{equation}
among the partial derivatives $f_j$'s of $f$ with coefficients $a_j$ in $S_q$, the vector space of  homogeneous polynomials of degree $q$. In this paper we assume $mdr(f)>0$, which is equivalent to asking $D$ not to be a cone over a hypersurface in $\PP^{n-1}$.

When $n=2$, then $D$ is a plane curve with isolated singularities and the degree of its singular subscheme $\Sigma$ is the global Tjurina number 
$$\tau(D)= \sum_{p \in \Sigma}\tau(D,p),$$
where $\tau(D,p)$ is the Tjurina number of the isolated plane curve singularity $(D,p)$, see for instance \cite{CD}. A. A. du Plessis and C.T.C. Wall have proved the following, see \cite{duPCTC}.

\begin{thm}
\label{thmI1}
In the class of reduced plane curves with a fixed degree $e_1=mdr(f)$  satisfying
$2e_1 \leq d-1$, the free curve $D$ (if it exists) has a singular locus $\Sigma$ of maximal degree. More precisely, for a curve $D$ with fixed degree $e_1$, one has
$$\deg \Sigma \leq s_1^2-s_2,$$
where $e_2=d-1-e_1$, and  $s_1=e_1+e_2$, $s_2=e_1e_2$  are the  elementary symmetric functions in $e_1,e_2$. Moreover,  the equality holds if and only if the curve $D$ is free.
\end{thm}

This result is also discussed in \cite{Dmax}, where a similar fact is proved for the nearly free plane curves.

In this note we investigate to what extent such  results hold for $n=3$, i. e.  for surfaces in $\PP^3$. 
The basic properties of a free surface $D$ are reviewed in the next section, for now we just say that freeness of $D$ is the same as asking the 3-fold singularity given by the cone over $D$ to be a free divisor germ in $(\C^4,0)$ in the sense of K. Saito, who introduced this important notion of free divisor in \cite{KS}.
Since free surfaces have a 1-dimensional singular locus, we restrict our attention to such surfaces, i.e. from now on we suppose 
$$\dim \Sigma=1.$$
Instead of looking only at the minimal degree $e_1$ of a Jacobian relation, we have to consider a pair of degrees $e_1 \leq e_2$, corresponding to the lowest degrees of two independent generators $\rho_1$ and $\rho_2$ of the module $AR(f)$ of all the Jacobian syzygies of $f$. The behaviour of these two generators (i.e. how independent they are) is used to define a class of surfaces, called tame surfaces, see Definition \ref{def3}. In particular, the free surfaces and most (if not all, see Question \ref{q1})  of the nearly free surfaces defined in \cite{DStFS}
are shown to be tame, see Proposition \ref{prop1}. The main results of this note are the following.

\begin{thm}
\label{thmI}
Let $D:f=0$ be a reduced  surface in $\PP^3$, tame with respect to the pair of syzygies $\rho_1,\rho_2$ of degrees $e_1 = \deg \rho_1 \leq e_2= \deg \rho_2$. Set $e_3=d-1-e_1-e_2$ and assume $e_3 \geq e_2$.
Then one has the following.

\begin{enumerate}

\item For any $k < e_3+d-1$, one has
$$\dim M(f)_k={k+3 \choose 3}-4{k-d+4 \choose 3}+\sum_{j=1,2}{ k-d-e_j+4 \choose 3}.$$
In particular, these dimensions are independent of $f$ once the degrees $e_1$ and $e_2$ are fixed.

\item For any $k \geq e_3+d-1$, one has
$$\dim M(f)_k \leq {k+3 \choose 3}-4{k-d+4 \choose 3}+\sum_{j=1,3}{ k-d-e_j+4 \choose 3},$$
and the equality holds for any $k \geq e_3+d-1$ if and only if it holds for $k = e_3+d-1$. Moreover, these equalities hold  exactly when $D:f=0$ is a free surface with exponents $e_1,e_2,e_3$.

\end{enumerate}

\end{thm}

\begin{thm}
\label{thmI2}
In the class of tame surfaces $D$ with fixed degrees $e_1,e_2$  satisfying
$e_1+2e_2 \leq d-1$, the free surface  (if it exists) has a singular locus $\Sigma$ of maximal degree. More precisely, for such a tame surface $D$ with fixed degrees $e_1,e_2$, one has
$$\deg \Sigma \leq s_1^2-s_2,$$
where $e_3=d-1-e_1-e_2$, and  $s_1=\sum_{j=1,3}e_j$, $s_2=\sum_{i <j}e_ie_j$  are the first two elementary symmetric functions in $e_1,e_2,e_3$. Moreover,  the equality holds if and only if the surface $D$ is free.
\end{thm}

The tameness assumption is necessary in both Theorems above, as follows from Example \ref{exNT}. 
Some simple characterizations of  the nearly free tame surfaces are also given, see Proposition \ref{propNF} and Corollary \ref{cor3}, in perfect analogy to the case of nearly free curves treated in \cite{Dmax}.

\medskip

We thank Aron Simis for useful discussions and Gabriel Sticlaru for finding the interesting surfaces described in Examples \ref{extame}, \ref{exNT} and \ref{exNF}.

\section{Free, nearly free and tame surfaces in $\PP^3$} 

Let $f$ be a homogeneous polynomial of degree $d$ in the polynomial ring $S=\C[x,y,z,w]$ and denote by $f_x,f_y,f_z, f_w$ the corresponding partial derivatives.
Let $D$ be the surface in $\PP^3$ defined by $f=0$ and assume that $D$ is reduced and not a cone over a plane curve. We denote by $J_f$ the Jacobian ideal of $f$, i.e. the homogeneous ideal in $S$ spanned by $f_x,f_y,f_z,f_w$, and  by $M(f)=S/J_f$ the corresponding graded ring, called the Jacobian (or Milnor) algebra of $f$.

 Consider the graded $S-$submodule $AR(f) \subset S^{4}$ of {\it all relations} involving the derivatives of $f$, namely
$$\rho=(\alpha, \beta,\gamma, \delta)\in AR(f)_q$$
if and only if  $\alpha f_x+ \beta f_y +\gamma f_z+ \delta f_w=0$ and $\alpha, \beta,\gamma, \delta$ are in $S_q$. We set $ar(f)_k=\dim AR(f)_k$ and    $m(f)_k=\dim M(f)_k$ for any integer $k$.

\begin{definition}
\label{def1}

 The surface $D:f=0$ is a { free divisor} if the following  equivalent conditions hold.

\begin{enumerate}

\item $M(f)$ is a Cohen-Macaulay $S$-module, i.e. $\depth M(f) =\dim M(f)=2$.

\item $H^0_{\bf m}(M(f))=H^1_{\bf m}(M(f))=0,$
 with ${\bf m}=(x,y,z,w)$ the maximal homogeneous ideal  in $S=\C[x,y,z,w]$.

\item The minimal resolution of the Milnor algebra $M(f)$ has the following  form
$$0 \to \oplus_{j=1,3} S(-d_j-d+1)  \to S^4(-d+1) \xrightarrow{(f_x,f_y,f_z,f_w)}  S$$
for some positive integers $d_1 \leq d_2 \leq d_3$.
\item The graded $S$-module $AR(f)$ is free of rank 3, i.e. there is an isomorphism 
$$AR(f)=S(-d_1) \oplus S(-d_2)\oplus S(-d_3) $$
for some positive integers $d_1 \leq d_2 \leq d_3$.
\item The coherent sheaf $\F$ on $\PP^3$ associated to the graded $S$-module $AR(f)$ splits as a direct sum of line bundles, i.e.
$$\F=\OO(-d_1) \oplus \OO(-d_2)\oplus \OO(-d_3) $$
for some positive integers $d_1 \leq d_2 \leq d_3$.
\end{enumerate}

\end{definition}
When $D$ is a free divisor, the integers $d_1 \leq d_2 \leq d_3$ are called the {exponents} of $D$.  They satisfy the relation $d_1+d_2+d_3=d-1 $ and the coefficients of the Hilbert polynomial $P(M(f))(k)=ak+b$ of the Milnor algebra $M(f)$ are given by
\begin{equation}
\label{free1}
a=s_1^2-s_2    \text{ and }  b=2a-s_1^3+\frac{3}{2}s_1s_2-\frac{1}{2}s_3,
\end{equation}
where $s_1=\sum_{j=1,3}d_j$, $s_2=\sum_{i <j}d_id_j$ and $s_3=d_1d_2d_3$ are the elementary symmetric functions in the exponents, see  \cite{DStFS}.

\begin{definition} 
\label{def2}

The surface $D:f=0$ is a { nearly free divisor} 
if the following  equivalent conditions hold.

\begin{enumerate}

\item The Milnor algebra $M(f)$ has a minimal resolution of the form
$$0 \to S(-d-d_3) \to \oplus_{j=1,4}S(-d-d_j+1)  \to S^4(-d+1) \xrightarrow{(f_x,f_y,f_z,f_w)}  S$$
for some integers $1 \leq d_1 \leq d_2\leq d_3=d_4$, called the exponents of $D$.

\item There are 4 syzygies $\rho_1$, $\rho_2$, $\rho_3$, $\rho_4$ of degrees $d_1 \leq d_2  \leq d_3=d_4=d-(d_1+d_2)$ which form a minimal system of generators for the first syzygies module $AR(f)$.

\end{enumerate}

\end{definition}

In down-to-earth terms, this definition says that the module $AR(f)$ is not free of rank 3, but it has 4 generators $\rho_1$, $\rho_2$, $\rho_3$ and $\rho_4$ of degree respectively $d_1$, $d_2$, $d_3$ and $d_3$ and the second order syzygy module is spanned by a unique relation
\begin{equation}
\label{SOS}
R: a_1\rho_1+a_2\rho_2+a_3\rho_3+a_4\rho_4=0,
\end{equation}
where $a_1,a_2,a_3,a_4$ are homogeneous polynomials in $S$ of degrees $d_3-d_1+1, d_3-d_2+1,1,1$ respectively.

If $D:f=0$ is nearly free, then the exponents $d_1 \leq d_2\leq d_3$ ($d_4$ is omitted since it coincides to $d_3$) determine the Hilbert polynomial $P(M(f))(k)=ak+b$ as follows.
Define $d'_1=d_1$, $d'_2=d_2$, $d'_3=d_3-1$ and let the integers $a'$ and $b'$ be computed using the formulas in \eqref{free1}, i.e. as if $a',b'$ were the coefficients of the Hilbert polynomial corresponding to a free surface $D'$ with exponents $d'_1,d'_2,d'_3$. Then one has the formulas 
\begin{equation}
\label{nfree1}
 a=a'-1 \text{   and   } b=b'+d+d_3-3,
\end{equation}
see \cite{DStFS}.
For both a free and a nearly free surface $D:f=0$, it is clear that $mdr(f)=d_1$.

 To an element $ \rho=(\alpha, \beta,\gamma, \delta) \in S^4$, we can associate the differential 1-form
$$\omega(\rho)=\alpha dx+ \beta dy+ \gamma dz+ \delta dw \in \Omega^1,$$
and consider the Koszul complex of $\alpha, \beta,\gamma, \delta$ in $S$ given by
$$K^*(\alpha, \beta,\gamma, \delta):  \   \   \  0 \to \Omega^0 \to \Omega^1 \to \Omega^2 \to \Omega^3 \to \Omega^4 \to  0,$$
where $ \Omega^k$ denotes the $S$-module of global algebraic differential $k$-forms on $\C^4$ and the morphisms are given by the wedge product by $ \omega(\rho)$. If we assume that $\alpha, \beta,\gamma, \delta$ do not have any common factor of degree $>0$, then
the grade of the ideal $I=(\alpha, \beta,\gamma, \delta)$ (i.e. the maximal length of a regular sequence contained in $I$), which  is equal to the codimension of $I$, is clearly at least 2. This implies the vanishing of the first cohomology of the Koszul complex $K^*(\alpha, \beta,\gamma, \delta)$, see for instance Thm. A.2.48 in \cite{Eis},
i.e. the following sequence
\begin{equation}
\label{es1}
 0 \to \Omega^0 \xrightarrow{\omega(\rho) }  \Omega^1 \xrightarrow{\omega(\rho) } \Omega^2
\end{equation}
is exact. This applies in particular for a surface $D:f=0$ if we choose $\rho=\rho_1$ to be a nonzero syzygy of minimal degree, say $e_1$, in $AR(f)$. Let now $\rho_2$ be a homogeneous representative in $AR(f)$ of a nonzero homogeneous element of minimal degree, say $e_2$ with $e_2 \geq e_1$, in the quotient module $AR(f)/(S\rho_1)$. Note that both $\rho_1$ and $\rho_2$ are {\it primitive syzygies}, i.e. they are not nonconstant multiples of lower degree syzygies.
Let 
$\Omega^1_0=\{\omega(\rho)  \ \ : \ \ \rho \in AR(f) \},$
and consider the sequence of graded $S$-modules
\begin{equation}
\label{es2}
  S(-e_1-1) \oplus S(-e_2-1)   \xrightarrow{u=(\omega(\rho_1), \omega(\rho_2) )}  \Omega^1_0 \xrightarrow{v=\omega(\rho_1) \wedge \omega(\rho_2) } \Omega^3(e_1+e_2+2),
\end{equation}
where the first morphism is $u:(a,b) \mapsto a\omega(\rho_1)+b\omega(\rho_2)$ and the second one $v$ is induced by the wedge product by $\omega(\rho_1) \wedge \omega(\rho_2)$.

\begin{definition} 
\label{def3}

The surface $D$ is tame with respect to the syzygies $\rho_1$ and $\rho_2$ if the sequence \eqref{es2} is exact.
The surface $D$ is tame if there is a pair of syzygies $\rho_1$ and $\rho_2$ as above, such that the surface is tame with respect to $\rho_1$ and $\rho_2$.

\end{definition}

Note that a 1-form 
$\omega(\rho)$ is in the kernel of $v$ if and only if  the matrix $M(\rho_1, \rho_2,\rho)$ with 3 rows, corresponding to the 4 components of $\rho_1, \rho_2$ and $\rho$ has rank 2 over the field of fractions $K=\C(x,y,z,w)$. The exactness of the sequence \eqref{es1} shows that this happens if and only if 
\begin{equation}
\label{eq10}
c \cdot \rho= c_1 \cdot \rho_1+c_2 \cdot \rho_2,
\end{equation}
for some polynomials $c,c_1,c_2$ in $S$ with $c \ne 0$.
This implies the following.

\begin{lem}
\label{lem1} (i)  The surface $D:f=0$ is not tame with respect to $\rho_1$ and $\rho_2$ exactly when there is a primitive syzygy $\rho$ satisfying the equality \eqref{eq10} with $\deg c>0$ and $c_1$ and $c_2$ relatively prime polynomials.

\noindent (ii) Suppose there is a closed Zariski subset $B \subset \C^4$ such that the matrix $M(\rho_1, \rho_2)$ with 2 rows, corresponding to the 4 components of $\rho_1$ and $\rho_2$ has rank 2 for any point $p=(x,y,z,w) \notin B$. Then the corresponding sequence \eqref{es2} is exact.

\end{lem}

The assumption in (ii) is equivalent to asking the six $2 \times 2$-minors of the matrix $M(\rho_1, \rho_2)$ not to have a common divisor in $S$.

\proof The first claim (i) is clear. To prove (ii),
assume that the sequence \eqref{es2} is not exact.
Note that the equality \eqref{eq10} implies that the rank of the matrix $M(\rho_1, \rho_2)$ is 1 on the zero set of $c$, which has codimension one since $c$ is not a constant, hence we get a contradiction.

\endproof

\begin{prop}
\label{prop1}

(i) Any free surface is tame.

\noindent (ii) A nearly free surface is tame if and only if the linear forms $a_3$ and $a_4$ which occur in the second order syzygy \eqref{SOS} are linearly independent.

\end{prop}

\proof (i) Indeed, if we choose $\rho_1$ and $\rho_2$ as the first two elements in the basis of $AR(f)$ described in Definition \ref{def1} (4), it is clear that a syzygy $\rho$ satisfies \eqref{eq10} if and only if it belongs to the image of $u$, i.e. it is a linear combination of 
$\rho_1$ and $\rho_2$ with coefficients in $S$.

(ii) If the linear forms $a_3$ and $a_4$ which occur in the second order syzygy \eqref{SOS} are linearly dependent, we may assume $a_4=0$ by choosing $\rho_3$ and $\rho_4$ appropriatedly.
Then $\rho=\rho_3$ satisfies $a_3\rho= -a_1\rho_1 -a_2 \rho_2$, which shows that $D:f=0$ is not tame with respect to $\rho_1$ and $\rho_2$. If $d_3>d_2$, then the choice for $\rho_1$ and $\rho_2$ is essentialy unique, and we are done. If $d_3=d_2>d_1$ or $d_3=d_2=d_1$, then the proof can be easily adapted by the reader, since the choices for $\rho_1$ and $\rho_2$ can be listed.

Suppose now that the linear forms $a_3$ and $a_4$ are linearly independent. Choose $\rho_1$ and $\rho_2$ as the first two elements in the system of minimal generators of $\rho_j$'s, with $j=1,...,4$ of $AR(f)$ described in Definition \ref{def2} (2) and write
$\rho$, an element in the kernel of $v$, as 
 $\rho= \sum_{i=1,4}b_i\rho_i$. Then Lemma \ref{lem1} (ii) implies the existence of polynomials $c,c_1,c_2,a$ such that $c_1, c_2$ are relatively prime and 
\begin{equation}
\label{SOS1}
aa_1=cb_1-c_1, \\ aa_2=cb_2-c_2, \\ aa_3=cb_3 \text {  and  } aa_4=cb_4.
\end{equation}
The first two equalities imply that $a$ and $c$ are relatively prime (since $c_1, c_2$ are relatively prime). Then the last two equalities imply that $c$ divides both $a_3$ and $a_4$, which is possible only if $c$ is a constant. Hence $D:f=0$ is  tame with respect to $\rho_1$ and $\rho_2$.
\endproof

All the examples of nearly free surfaces described in \cite{DStFS} are tame, because the linear forms $a_3$ and $a_4$ which occur in the second order syzygy \eqref{SOS} are linearly independent. This can be seen also by applying to each example Lemma \ref{lem1} (ii).
It is natural to ask the following.

\begin{question}
\label{q1}
Does there exist a nearly free surface in $\PP^3$ which is not tame?
\end{question}
Note that for a nearly free plane curve, the corresponding two linear forms occuring in the second syzygy similar to \eqref{SOS} are always linearly independent, see Remark 5.2 in \cite{DStFS}.

\begin{ex}
\label{extame} 

\noindent (i) We give now an example of a tame surface which is neither free nor nearly free, and for which the choice of the syzygies $\rho_1$ and $\rho_2$ is more subtle.
Let $D:f=(xw+y^2)^2+y^2z^2=0$.
Then one has the following generating syzygies:

\medskip

\noindent $\syz[1]: (x)f_x+(0)f_y+(0)f_z+(-w)f_w=0,$ 

\medskip
\noindent $\syz[2]: (2y^2)f_x+(-yw)f_y+(zw)f_z+(0)f_w=0,$

\medskip
\noindent $\syz[3]: (0)f_x+(-yz)f_y+(2y^2+z^2+2xw)f_z+(0)f_w=0,$

\medskip
\noindent $\syz[4]: (0)f_x+(xy)f_y+(-xz)f_z+(-2y^2)f_w=0$

\medskip
\noindent and some higher degrees ones. If we choose $\rho_1=\syz[1]$ and $\rho_2=\syz[2]$,  the sequence \eqref{es2} is not exact. Indeed, one has the following relation
$$w \rho=2y^2 \rho_1 -x\rho_2,$$
where $\rho=\syz[4]$, compare with the equality \eqref{eq10}. On the other hand, if we choose
$\rho_1=\syz[1]$ and $\rho_2=\syz[3]$, then it is easy to see that the rank of matrix $M(\rho_1, \rho_2)$ is 2 outside a subset of codimension 2. Then the corresponding sequence \eqref{es2} is exact by Lemma \ref{lem1}. In conclusion, the surface $D$ is tame.

\medskip

\noindent (ii) Here we give an example of a surface which is not tame.
Let $D: f=(x^2+y^2+zw)^4+y^4z^4=0$. This surface satisfies $H^0_{\bf m}(M(f))=0$ and has the following generators for $AR(f)$

\medskip
\noindent $\syz[1]: (-z)f_x+(0)f_y+(0)f_z+(2x)f_w=0,$

\medskip
\noindent $\syz[2]: (0)f_x+(-yz)f_y+(z^2)f_z+(2y^2-zw)f_w=0,$

\medskip
\noindent $\syz[3]: (-2y^2+zw)f_x+(2xy)f_y+(-2xz)f_z+(0)f_w=0,$

\medskip
\noindent plus some higher degree ones. The relation
$$(2y^2-zw)\syz[1]-2x \syz[2]=z\syz[3]$$
shows that the only possible choices, i.e. $\rho_1=\syz[1]$ and $\rho_2$ a linear combination of $\syz[2]$ and $\syz[3]$, cannot produce an exact sequence due to an equality as in \eqref{eq10}.

\medskip

\noindent (iii) Finally we describe a plane arrangement which is not tame. Consider the arrangement
$$D:  f=w(x+y)(y+z)(x+z)(y-2z)(x+2y+3z)(11x+7y+5z+3w).$$
Then $AR(f)$ has as generating syzygies one syzygy $\rho_1$ of degree 2, two syzygies $\rho_2$ and $\rho_3$ of degree 3, and some other higher degree generators. The matrix $M(\rho_1, \rho_2,\rho_3)$ has rank 2 over the field $K=\C(x,y,z,w)$, and this implies that $D$ is not tame.

\end{ex}

A natural question is the following.

\begin{question}
\label{q2}
Can the tameness of  a plane arrangement in $\PP^3$ be characterized combinatorially?
\end{question}
A positive answer to this question would have implications for the Terao's conjecture (see \cite{Yo} for a discussion of this conjecture), similar to those for line arrangements discussed in \cite{Dmax}.
This comes from the fact that $\deg \Sigma$ is known to be determined by the combinatorics, see 
\cite{HS} and \cite{WY}.

\section{Bourbaki ideal of the syzygy module} 

For a reduced surface $D:f=0$ in $\PP^3$, we choose $\rho=\rho_1$ to be a nonzero syzygy of minimal degree, say $e_1$, and  $\rho_2$  a homogeneous representative in $AR(f)$ of a nonzero homogeneous element of minimal degree, say $e_2$ with $e_2 \geq e_1$, in the quotient module $AR(f)/(S\rho_1)$.

Let $X=\nabla f$ be  gradient vector field of $f$ on $\C^4$ and denote by $\iota_X:\Omega^k \to \Omega^{k-1}$ the interior product given by the contraction of a differential fom by the vector field $X$.
The (homology) complex
\begin{equation}
\label{KC2}
0 \to \Omega^4 \xrightarrow{ \iota_X}   \Omega^3 \xrightarrow{ \iota_X} \Omega^2 \xrightarrow{ \iota_X} \Omega^1  \xrightarrow{ \iota_X} \Omega^0 \to 0,
\end{equation}
is nothing else but the Koszul complex of the partial derivatives $f_x,f_y,f_z,f_w$ of $f$. Since $D$ is reduced, it has a singular set of codimension at most 2 in $\PP^3$, and we get, exactly as we have obtained the exact sequence \ref{es1}, an exact sequence
\begin{equation}
\label{es3}
0 \to \Omega^4 \xrightarrow{ \iota_X}   \Omega^3 \xrightarrow{ \iota_X} \Omega^2.
\end{equation}
Now note that a 1-form $\omega$ is in $\Omega^1_0$, i.e. it comes from a syzygy in $AR(f)$, if and only if $ \iota_X(\omega)=0$. Since $ \iota_X$ is an anti-derivation, i.e. it satisfies a graded Leibnitz rule, it follows that the image of the morphism $v$ from \eqref{es2} is contained in
$$\Omega^3_0=\{\omega \in \Omega^3 \ \ : \ \  \iota_X(\omega)=0\}.$$
And the exact sequence \eqref{es3} gives an isomorphism 
$$S \cong \Omega^4(4) \cong \Omega^3_0(d+2).$$
The above proves the following result.

\begin{thm}
\label{thm1}
Let $D:f=0$ be a reduced  surface in $\PP^3$, tame with respect to the pair $\rho_1,\rho_2$ having degrees $e_1 = \deg \rho_1 \leq e_2= \deg \rho_2$. Set $e_3=d-1-e_1-e_2$.
Then one has the following exact sequence of graded $S$-modules
$$ 0 \to  S(-e_1) \oplus S(-e_2)   \xrightarrow{u'} AR(f) \xrightarrow{v'} S(-e_3),$$
where $u'(a,b)=a\rho_1+b\rho_2$ and $v'(\rho)=h$, with $h$ the unique polynomial such that
$$\omega(\rho_1) \wedge \omega(\rho_2)  \wedge \omega(\rho) =h\iota_X(dx\wedge dy \wedge dz \wedge dw).$$
\end{thm}

If we denote by $B(f) \subset S$ the image of the morphism $v'$, it follows that $B(f)$ is a Bourbaki ideal for the syzygy module $AR(f)$, see \cite{Bou}, Chapitre 7, \S 4, Thm. 6, as well as section 3 in \cite{SUV}.

\begin{cor}
\label{corBou}
With the above notation, the minimal number of generators of the Bourbaki ideal $B(f)$ is the same as the minimal number of syzygies that one must add to $\rho_1$ and $\rho_2$ in order to get a generating set of the syzygy module $AR(f)$.
\end{cor}

\begin{rk}
\label{rkSaito}
If we denote by $M(E,\rho_1,\rho_2,\rho)$ the $4 \times 4$ matrix having as the first row the components of the Euler vector field, i.e. $(x,y,z,w)$, and as the next rows the components of $\rho_1,\rho_2$ and $\rho$ respectively, then the equality $v'(\rho)=h$ is equivalent to
$$\det M(E,\rho_1,\rho_2,\rho)=h\cdot f.$$
\end{rk}

By the choice of the syzygies $\rho_1,\rho_2$, it is clear that $(\im u')_k=AR(f)_k$, for any $k\leq e_2-1$.  Then Theorem \ref{thm1} implies the following consequence, which is essentially the same as Lemma 4.12 in \cite{DStFS}, which in turn is an extension of the corresponding result for free curves in Lemma 1.1 in \cite{ST}. This result also shed a new light on Saito's criterion of freeness, see \cite{KS}, \cite{Yo}, \cite{DStFS}, which is used to get the freeness of $D$ in the result below.

\begin{cor}
\label{cor1} 

Let $p$ be the smallest integer such that $(\im u')_p \ne AR(f)_p$, i.e. $p$ is the first degree where some new generator for $AR(f)$ should be added besides $\rho_1$ and $\rho_2$. Assume $e_3=d-1-e_1-e_2 \geq e_2$.
Then $p \geq e_3$ and the following conditions are equivalent.

\begin{enumerate}

\item $p = e_3$;

\item $D:f=0$ is a free surface with exponents $e_1,e_2,e_3$;

\item $$ar(f)_{e_3}> { e_3-e_1+3 \choose 3} +{ e_3-e_2+3 \choose 3};$$

\item $$ar(f)_{e_3}={ e_3-e_1+3 \choose 3} +{ e_3-e_2+3 \choose 3}+1.$$

\item For any $k \geq e_3$ one has $$ar(f)_{k}\geq { k-e_1+3 \choose 3} +{ k-e_2+3 \choose 3}+  { k-e_3+3 \choose 3};$$

\item For any $k \geq e_3$ one has $$ar(f)_{k}= { k-e_1+3 \choose 3} +{ k-e_2+3 \choose 3}+  { k-e_3+3 \choose 3}.$$

\end{enumerate}

\end{cor}

\section{Proof of Theorems \ref{thmI} and  \ref{thmI2}} 

We have an obvious exact sequence
$$0 \to AR(f)_{k-d+1} \to S^4_{k-d+1} \to S_k \to M(f)_k \to 0.$$
This implies
\begin{equation}
\label{mf1}
m(f)_k={k+3 \choose 3}-4{k-d+4 \choose 3}+ar(f)_{k-d+1},
\end{equation}
for $k \geq d-4$.

 Theorem \ref{thm1} and Corollary \ref{cor1} clearly imply Theorem \ref{thmI}.

\begin{ex}
\label{exNT} For the surface $D:f=0$ in Example \ref{extame}, we have $\deg \Sigma=36$, while 
$s_1^2-s_2=35.$ Indeed, in this case $e_1=1$, $e_2=2$ and $e_3=d-1-e_1-e_2=4$.
Hence the condition of tameness is  necessary to have the result in Corollary \ref{thmI2}.

\end{ex}

The proof of Theorem   \ref{thmI2} is more involved, even if in view of Theorem \ref{thmI},  the only point that needs explanation is the fact that the equality $s_1^2-s_2=\deg \Sigma$ implies that $D:f=0$ is  free.
Consider the following exact sequence of graded $S$-modules
\begin{equation}
\label{BES}
0 \to  S(-e_1) \oplus S(-e_2)   \xrightarrow{u'} AR(f) \xrightarrow{v'} B(f)(-e_3)\to 0,
\end{equation}
coming from Theorem \ref{thm1}, where $B(f)$ is the Bourbaki ideal of $AR(f)$. Using this and the obvious exact sequence
$$0 \to B(f)(-e_3) \to S(-e_3) \to (S/B(f))(-e_3)\to 0,$$
we get
$$ar(f)_{k}= { k-e_1+3 \choose 3} +{ k-e_2+3 \choose 3}+  { k-e_3+3 \choose 3} -\dim (S/B(f))_{k-e_3}.$$
Combining this equality with the relation \eqref{mf1} written for $k$ replaced by $k+d-1$,
we get, for $k$ large enough,
\begin{equation}
\label{mf2}
\dim (S/B(f))_{k-e_3}=Q(k)-m(f)_{k+d-1},
\end{equation}
where the polynomial $Q(k)$ is given by
\begin{equation}
\label{mf3}
Q(k)=\sum_{j=1,3}{ k-e_j+3 \choose 3}-4{ k+3 \choose 3}+{ k+d+2 \choose 3}.
\end{equation}
If we expand binomial coefficients, we get that $Q(k)$ is a polynomial in $k$ of degree (at most) one and the coefficient of $k$ is given by $s_1^2-s_2$. It follows that the equality $s_1^2-s_2=\deg \Sigma$ and the equality \eqref{mf2} imply that the subscheme $Y$ of $\PP^3$ defined by the ideal $B(f)$ is 0-dimensional or empty. On the other hand, Lemma 3 in \cite{BR} shows that the existence of the exact sequence \eqref{BES} with $AR(f)$ a reflexive $S$-module and $\codim Y\geq 2$ imply that either $Y$ has pure codimension two, or $Y$ is empty. See also the discussion on p. 145 in \cite{BEG} and note that $AR(f)$ is a second syzygy module by definition. Hence the only possibility is that $Y$ is empty, i.e. $B(f)$ is an $\bf m$-primary ideal. Then the exact sequence \eqref{BES} translates into the following exact sequence of sheaves on $\PP^3$
$$0 \to \OO(-e_1) \oplus \OO(-e_2) \to \F \to \OO(-e_3) \to 0,$$
where $\F$ is the sheaf associated to the graded $S$-module $AR(f)$.
But this implies that $ \F=\OO(-e_1) \oplus \OO(-e_2) \oplus \OO(-e_3)$, i.e. $D:f=0$ is free with exponents $e_1,e_2,e_3$, since the group of extensions of $\OO(-e_3)$ by $\OO(-e_1) \oplus \OO(-e_2)$ is
$$\Ext ^1(\OO(-e_3),\OO(-e_1) \oplus \OO(-e_2))=\Ext ^1(\OO,\OO(e_3-e_1) \oplus \OO(e_3-e_2))=$$
$$=H^1(\PP^3,\OO(e_3-e_1) \oplus \OO(e_3-e_2))=0.$$
This completes the proof of Theorem   \ref{thmI2}. 

\section{Two characterizations of nearly free tame surfaces} 

\begin{prop}
\label{propNF}

Let $D:f=0$ be a reduced  surface in $\PP^3$, tame with respect to the pair of syzygies $\rho_1,\rho_2$ of degrees $e_1 = \deg \rho_1 \leq e_2= \deg \rho_2$. Set $e_3=d-1-e_1-e_2$ and assume $e_3 \geq e_2$.
Then $D$ is nearly free with exponents $e_1,e_2,e_3+1$ if and only if $\dim B(f) _{1}=2$.

\end{prop}
\proof If $D$ is nearly free, then the equality $\dim B(f) _{1} =2$ follows from Proposition \ref{prop1}, (ii), since
the image under $v'$ of the syzygies $\rho_3$ and $\rho_4$ are, up to a constant factor, $a_4$ and $-a_3$, in the notation from the proof of Proposition \ref{prop1}.
Conversely, choose $\rho_3$ and $\rho_4$ such that their images $a'_3$ and $a'_4$ under $v'$ span $B(f) _{1}$. It follows that the subscheme $Y$ defined by the ideal $B(f)$ is contained in the projective line $L$ in $\PP^3$ defined by $a'_3=a'_4=0$. As in the proof above of Theorem \ref{thmI2}, using Lemma 3 in \cite{BR}, we see that there are only the following two cases to discuss.

\noindent {\bf Case 1.} $Y$ is empty, and this corresponds as we have seen in Theorem \ref{thmI2} to the case $D$ is a free divisor. But for a free divisor one has $\dim B(f)_1=4$, as the image of $v'$ is the whole ring $S$.
Hence this case is impossible.

\noindent {\bf Case 2.} $Y$ has pure codimension 2. It follows that $Y=L$, and hence $B(f)$ is the ideal spanned by the independent linear forms $a'_3$ and $a'_4$. This implies using Theorem \ref{thm1} that the four syzygies 
$\rho_1$, $\rho_2$, $\rho_3$ and $\rho_4$ span the syzygy module $AR(f)$, and hence $D:f=0$ is a nearly free surface.

\endproof

\begin{rk}
\label{rkNF}
(i) To check that a given surface is tame with respect to the pair of syzygies $\rho_1,\rho_2$ without having detailed information on the syzygy module $AR(f)$
can be done using Lemma \ref{lem1} (ii).

\noindent (ii) In Proposition \ref{propNF}, the condition $\dim B(f) _{1} =2$ can be replaced by the condition
$\dim B(f) _{1} \geq 2$ and $D$ is not free. This is the perfect analog of the characterization of nearly free curves in Theorem 4.1 (ii) in \cite{Dmax}, except that for surfaces we need the extra tameness condition.
\end{rk}

\begin{cor}
\label{cor3}
Let $D:f=0$ be a reduced  surface in $\PP^3$, tame with respect to the pair of syzygies $\rho_1,\rho_2$ of degrees $e_1 = \deg \rho_1 \leq e_2= \deg \rho_2$. Set $e_3=d-1-e_1-e_2$ and assume $e_3 \geq e_2$ and that $D$ is not free. Then
$$\deg \Sigma \leq s_1^2-s_2-1,$$
and the equality holds if the surface $D$ is nearly free  with exponents $e_1,e_2,e_3+1$.
Conversely, if this equality holds and if the syzygy module $AR(f)$ is spanned by $\rho_1,\rho_2$ and two other syzygies $\rho_3$ and $\rho_4$, then $D$ is a nearly free surface with exponents $e_1,e_2,e_3+1$.
\end{cor}

\proof
The only point that needs explanation is the last claim. The formula \eqref{mf2} implies that $\deg Y$, which is the leading coefficient of the Hilbert polynomial $H(S/B(f))$, is one. On the other hand, we know as above that $Y$ has pure codimension 2, since it is non empty. Under our assumption, $Y$ is in fact a complete intersection, since the ideal $B(f)$ is spanned by $g_j=v'(\rho_j)$ for $j=3,4$. The degree of the complete intersection $Y$ is the product $\deg g_3 \cdot \deg g_4$, hence both $g_3$ and $g_4$ are linear forms.
The claim then follows from Proposition \ref{propNF}.

\endproof

The following example shows that the equality $s_1^2-s_2-1=\deg \Sigma$ does not imply that $D:f=0$ is nearly free without the tameness assumption.
\begin{ex}
\label{exNF}
Consider the surface
$D:f =x^{12}+x^{11}y+ y^{11}z+w(x^{10}w+y^{10}z).$ This surface is not nearly free since $H^0_{\bf m}(M(f)) \ne 0$,
see \cite{DStFS}. The generators of $AR(f)$ have degrees $\geq 3$ and $\dim AR(f)_3=4$.
The corresponding Hilbert polynomial is
$H(M(f))(k) = 81k - 507$   for $k \geq 21$. If we set $e_1=e_2=3$ and $e_3=d-1-e_1-e_2=5$, we get
$$s_1^2-s_2-1=81= \deg \Sigma.$$
\end{ex}

\end{document}